\documentclass[final,onefignum,onetabnum]{article}

\usepackage{array}

\usepackage[caption=false]{subfig}
\usepackage{pgfplots}
\usepackage{cite}

\usepackage{algorithmic}

\usepackage{graphicx,epstopdf}

\usepackage{subfig}
\usepackage{amsmath}
\usepackage{amssymb}
\usepackage{graphicx}
\usepackage{dcolumn}
\usepackage{mathtools}
\usepackage{amsmath,amsfonts}
\usepackage{bm}
\usepackage{amsmath}
\usepackage{amssymb}
\usepackage{color}
\usepackage{float}
\usepackage{setspace}
\usepackage{tabularx}

\allowdisplaybreaks

\newcommand{\p}{\partial}
\newcommand{\f}{\frac}
\newcommand{\ds}{\displaystyle}
\newcommand{\E}{ {\mathbb{E}} }
\newcommand{\ola}{\overleftarrow}

\newcommand{\be}{\begin{equation}}
\newcommand{\ee}{\end{equation}}

\def\ba{\begin{array}}                \def\ea{\end{array}}
\def\bel{\begin{equation}\label}      \def\ee{\end{equation}}

\colorlet{texcscolor}{blue!50!black}
\colorlet{texemcolor}{red!70!black}
\colorlet{texpreamble}{red!70!black}
\colorlet{codebackground}{black!25!white!25}

\date{}

\title{A PDE-based Adaptive Kernel Method for Solving Optimal Filtering Problems}

\author{
Zezhong Zhang \thanks{ Department of Mathematics, Florida State University, Tallahassee, Florida.} 
\and Richard Archibald\thanks{ Computational Science and Mathematics Division, Oak Ridge National Laboratory, Oak Ridge, Tennessee,  \ ({\tt archibaldrk@ornl.gov}). }
\and Feng Bao\thanks{ Department of Mathematics, Florida State University, Tallahassee, Florida, \ ({\tt bao@math.fsu.edu}).} 
      }

\begin{document}
\maketitle

\begin{abstract}
In this paper, we introduce an adaptive kernel method for solving the optimal filtering problem. The computational framework that we adopt is the Bayesian filter, in which we recursively generate an optimal estimate for the state of a target stochastic dynamical system based on partial noisy observational data. The mathematical model that we use to formulate the propagation of the state dynamics is the Fokker-Planck equation, and we introduce an operator decomposition method to efficiently solve the Fokker-Planck equation. An adaptive kernel method is introduced to adaptively construct Gaussian kernels to approximate the probability distribution of the target state. Bayesian inference is applied to incorporate the observational data into the state model simulation. Numerical experiments have been carried out to validate the performance of our kernel method. 

\end{abstract}
\vspace{1em}

\vspace{1em}



\section{Introduction}
Data assimilation is an important topic in data science. It aims to optimally combine a mathematical model with observational data. The key mission in data assimilation is the optimal filtering problem, in which we try to find the best estimate for the state of a stochastic dynamic model. Such a stochastic dynamic model is typically in the form of a system of stochastic differential equations (SDEs), and we call it the ``state process''. In many practical situations, the true value of the state process is not available, and we can only use partial noisy observations to find the best estimate for the state. In the optimal filtering problem, the ``best estimate'' that we want to find is defined as the conditional expectation of the state conditioning on the observations. 

\vspace{0.3em}

When both the state dynamics and the observations are linear, the optimal filtering problem is a linear filtering problem, which can be analytically solved by the classic Kalman filter \cite{Kalman1961}. In the case of nonlinear filtering problem, one needs to derive an approximation for the conditional probability of the state -- instead of calculating the conditional expectation directly, and we call this conditional probability the ``filtering density''. There are two well-known nonlinear filtering methods, e.g., the Zakai's approach and the Bayesian filter. The Zakai's approach formulates the filtering density as the solution of a parabolic type stochastic partial differential equation \cite{zakai}, e.g., the Zakai equation, then we solve the  Zakai equation numerically to obtain an approximation for the filtering density \cite{Bao_Zakai, Zhang_Zakai, Gobet-Zakai, Bao_Half}. The Bayesian filter solves the nonlinear filtering problem in a recursive two-step procedure: a prediction step and an update step. In the prediction step, we predict the filtering density of the state before reception of the observational data. Then, when the observational data is available in the update step, we let the predicted filtering density be the prior and apply Bayesian inference to incorporate the observational information into the prediction and obtain an updated filtering density as the posterior in the Bayesian inference. Besides the Zakai's approach and the Bayesian filter, recently we developed a new approach to solve the optimal filtering problem, and we named this approach the backward SDE filter \cite{Bao_CMS, BSDE_filter, Bao_CiCP_2019}. The backward SDE filter is similar to Zakai's approach in the sense that it utilizes a system of differential equations (e.g., the backward SDEs) to analytically propagate the filtering density \cite{BCZ_2018, Bao_filter_jump}. At the same time, taking advantage of the SDE nature, the backward SDE filter is highly scalable, which can be implemented efficiently on parallel computers.
\vspace{0.3em}

Among the aforementioned approaches, the Bayesian filter is widely used for solving the nonlinear filtering problem in practice. The state-of-the-art Bayesian filter methods include the particle filter \cite{particle-filter, MCMC-PF} and the ensemble Kalman filter \cite{Evense_EnKF, EvensenBook}.   Both the particle filter and the ensemble Kalman filter use samples (particles) to create an empirical distribution to describe the predicted filtering density (i.e. the prior) in the prediction step. In the update step, the particle filter applies Bayesian inference to assign weights to particles and use weighted particles to represent the updated filtering density ( i.e. the posterior). On the other hand, the ensemble Kalman filter linearizes the observations and then adopts the Kalman update in the Kalman filter method to derive an updated filtering density. Since the Kalman filter is designed for linear problems, in the case that the optimal filtering problem is highly nonlinear, the ensemble Kalman filter does not provide accurate estimates for the target state \cite{DA-applications, vanLeeuwen}. The major drawback of the particle filter is the degeneracy issue \cite{Kang-PF, Sny-particle}. When the observational data lies on the tail of the predicted filtering density, only very few particles will receive high likelihood weights in the Bayesian inference procedure, which significantly reduces the effective particle size in the particle filter.  
\vspace{0.3em}

In addition to the methodology of using particles to propagate the state process, another approach that can transport the probability density forward in time is to solve the Fokker-Planck equation, which is a parabolic type particle differential equation (PDE). Although the PDE-based Fokker-Planck approach analytically formulates the  propagation of the state model, solving a PDE in high dimensional state space is computationally expensive, which makes PDE-based optimal filtering solvers difficult to be practically applied.  

\vspace{0.3em}

In this work, we will develop a novel kernel method to efficiently solve the Fokker-Planck equation, and the kernel approximated solution for the Fokker-Planck equation will be used as our estimate for the predicted filtering density in the prediction step. Then, we will adopt Bayesian inference to incorporate the observational data into the kernel approximated filtering density. 
\vspace{0.3em}

Kernel method recently attracted extensive attentions in machine learning and function approximation \cite{Kernel_learning, Multi-Kernel_Learning}. When solving the optimal filtering problem, the target function that we approximate with kernels is the filtering density, which is a probability density function (PDF). In many scenarios in the optimal filtering problem, the filtering density appears to be a bell-shaped function. This makes kernel method (especially with Gaussian type kernels) an effective way to construct approximations for the target filtering density \cite{Bao_Kernel_Arxiv, Bao_Kernel}. Since optimal filtering is often used to solve practical application problems in real time, efficiency of an optimal filtering method is essential. In this paper, we will introduce an operator decomposition method to decompose PDF propagation in the Fokker-Planck equation into a linear component and a nonlinear component. The linear component of propagation can be analytically derived, and the nonlinear component needs to be carried out numerically. Numerical solver for the nonlinear component in the Fokker-Planck equation will be formulated as an optimization problem, which aims to determine kernel parameters that describe the nonlinear propagation of the filtering density. To implement the optimization procedure efficiently, we will introduce a boosting algorithm \cite{NIPS1999_Boosting} to adaptively generate kernels to capture the main features of the state distribution. This allows us to use minimum amount of active kernels to characterize the filtering density, which will be used to estimate the target state. 

\vspace{0.3em}

The rest of this paper is organized as follows. In Section \ref{Prelim}, we introduce some preliminaries that we need to design the PDE-based adaptive kernel method for solving the optimal filtering problem. Then, we shall give a detailed description for our adaptive kernel method in Section \ref{Methodology}. Numerical examples that validate the effectiveness of the kernel approach in solving the optimal filtering problem and comparison experiments will be presented in Section \ref{Numerics}. Finally, summary and concluding remarks will be given in Section \ref{Conclusion}.


\section{Preliminaries}\label{Prelim}
In this section, we provide the preliminaries to formulate our adaptive kernel method for solving the optimal filtering problem. We shall first briefly introduce the optimal filtering problem. Then, we will discuss one of the most important optimal filtering approaches, i.e., the Bayesian filter, and we will describe the mathematical framework of the adaptive kernel method as a Bayesian filter type approach.

\subsection{The optimal filtering problem}
In the optimal filtering problem, we consider the following stochastic dynamical system in the form of a stochastic differential equation (SDE) in the probability space $(\Omega, \mathcal{F}, \mathbb{P})$
\begin{equation}\label{State}
dX_t = b(t, X_t) dt + \sigma_t dW_t,
\end{equation} 
where $b: \mathbb{R}^{+} \times \mathbb{R}^d \rightarrow \mathbb{R}^d$ is the drift coefficient, $\sigma: \mathbb{R}^+ \rightarrow \mathbb{R}^{d\times r}$ is the diffusion coefficient of the SDE, $W$ is a standard $r$-dimensional Brownian motion under $\mathbb{P}$, and the $\sigma_t dW_t$ term is a standard It\^o type stochastic integral, which brings additive noises to the dynamical model.  The $d$-dimensional stochastic process $X := \{X_t\}_{t \geq 0}$ is called the ``state process'', which represents the state of the dynamical model. In order to estimate the state of $X_t$ when the true value of $X_t$ is not available, we collect partial noisy observational data for $X_t$, denoted by $Y_t$, which is defined by
\begin{equation}\label{Observation}
Y_t = h(X_t) dt + dB_t,
\end{equation} 
where $h: \mathbb{R}^d \rightarrow \mathbb{R}^l$ is an observation function that measures the state of $X_t$ and $B$ is another Brownian motion independent of $W$ with covariance $R$ at any given time $t$. The stochastic process $Y$ is often called the ``observation process''. 

The goal of the optimal filtering problem is to find the best estimate for $\Psi(X_t)$ given the observational information $\mathcal{Y}_t$, where $\mathcal{Y}_t:=\sigma(Y_s, 0 \leq s \leq t)$ is the $\sigma$-algebra generated by the observation process $Y$, and $\Psi$ is a given test function. In mathematics, the best estimate for $\Psi(X_t)$ is defined by the ``optimal filter'', denoted by $\tilde{\Psi}(X_t)$, which is the conditional expectation of $\Psi(X_t)$, i.e. 
\begin{equation}\label{filter}
\tilde{\Psi}(X_t) : = \E[\Psi(X_t) | \mathcal{Y}_t]. 
\end{equation}

In this paper, we focus on the case that $f$ (in the state process) and/or $h$ (in the observation process) are nonlinear functions. The linear filtering problem is well-solved by the Kalman filter (except for the extremely high dimensional cases).
To solve the nonlinear optimal filtering problem, the standard approach aims to estimate the conditional probability of the state, i.e. $P(X_t | \mathcal{Y}_t)$, which is also called the ``filtering density''. Then, we can calculate the conditional expectation in Eq. \eqref{filter} through the integration formula
$$\E[\Psi(X_t) | \mathcal{Y}_t] = \int \Psi(x) P(x | \mathcal{Y}_t) dx. $$

In what follows, we will introduce the Bayesian filter, which provides a two-step procedure to estimate the filtering density  $P(X_t | \mathcal{Y}_t)$ recursively.

\subsection{The recursive Bayesian filter}

The Bayesian filter recursively estimates the target state $X_t$ on a sequence of discrete time instants $0 = t_1 < t_2 < \cdots, < t_n < \cdots$, and the Bayesian filter framework is composed of two steps: the \textit{prediction step} and the \textit{update step}. 
\vspace{0.5em}

\noindent \textbf{Prediction step.}
\vspace{0.5em}

Assume that the filtering density $p(X_{t_n} | \mathcal{Y}_{t_n})$ is given at time $t_n$. In the prediction step, we propagate the filtering density from time $t_n$ to time $t_{n+1}$ without usage of the new observational data $Y_{t_{n+1}}$, and we want to get the \textit{predicted filtering density}, i.e. $p(X_{t_{n+1}} | \mathcal{Y}_{t_n})$. 
\vspace{0.25em}

There are three major methods to achieve this goal:
\vspace{0.25em}

The first method is designed to find the predicted filtering density through the following Chapman-Kolmogorov formula 
$$p(X_{t_{n+1}} | \mathcal{Y}_{t_n}) = \int p(X_{t_{n+1}} | X_{t_n}) p(X_{t_n} | \mathcal{Y}_{t_n}) d X_{t_n},$$
where $p(X_{t_{n+1}} | X_{t_n})$ is the transition probability of the state equation \eqref{State} that transports the previous filtering density $p(X_{t_n} | \mathcal{Y}_{t_n})$ from $t_{n}$ to $t_{n+1}$. The above Chapman-Kolmogorov formula is often carried out by independent sample simulations, and it's the primary prediction technique in particle-based optimal filtering methods, such like the particle filter and the ensemble Kalman filter. As a result of the particle propagation of the filtering density, one may obtain an empirical representation for the predicted filtering density.

The second method utilizes the following (time-inverse) backward stochastic differential equation (BSDE) to generate the predicted filtering density:
\begin{equation*}
P_{t_{n+1}} = P_{t_n} - \int_{t_n}^{t_{n+1}} \sum_{i=1}^d \f{\p b_i}{\p x_i}(X_t) P_t dt - \int_{t_n}^{t_{n+1}} Q_t d\ola{W}_t, \qquad P_{t_n} = p(X_{t_n} | \mathcal{Y}_{t_n}),
\end{equation*}
where $X_t$ is the state process, and the $\int_{t_n}^{t_{n+1}} \cdot d\ola{W}_t $ is a backward It\^o integral, which is an It\^o type stochastic integral integrated backwards \cite{PP1994, Bao_Semi}. The solutions of the above BSDE is a pair $(P, Q)$, where $Q$ is the martingale representation of $P$ with respect to $W$ \cite{BSDE_finance}.
We refer to \cite{Bao_CMS, BCZ_2018, BSDE_filter, Bao_first} for more details of the BSDE method. 

The third method, which is also the method that we are going to discuss in this paper, describes the propagation of the filtering density through the following Fokker-Planck equation over the time interval $[t_{n}, t_{n+1}]$
\begin{equation}\label{Fokker-Planck}
\f{\p p(x, t)}{\p t} = - \sum_{i=1}^d \f{\p}{\p x_i}\big[ b_i(x, t) p(x, t)\big] + \sum_{i, j=1}^d \f{\p^2}{\p x_i \p x_j} \big[ D_{i, j} \ p(x, t)\big]
\end{equation}
with initial condition $p(x, t_{n}) = p(X_{t_n} = x | \mathcal{Y}_{t_n})$, where $b_i$ is the $i$-th component of the drift function $b$, and the matrix $D$ is defined by $D = \f{1}{2} \sigma \sigma^{\top}$. As a result, solution $p(x, t_{n+1})$ of the Fokker-Planck equation \eqref{Fokker-Planck} gives us the desired predicted filtering density $p(X_{t_{n+1}} = x | \mathcal{Y}_{t_n})$.

\vspace{0.5em}
\noindent \textbf{Update step.}
\vspace{0.5em}

With an approximation for the predicted filtering density (obtained through either one of the aforementioned method), the Bayesian filter updates the predicted filtering density to the (posterior) filtering density via the following Bayesian inference formula
\begin{equation}\label{Bayes}
p(X_{t_{n+1}}|\mathcal{Y}_{t_{n+1}}) = \f{p(X_{t_{n+1}}|\mathcal{Y}_{t_{n}}) \cdot p(Y_{t_{n+1}}|X_{t_{n+1}})}{p(Y_{t_{n+1}} | \mathcal{Y}_{t_n}) },
\end{equation}
where 
\begin{equation}\label{Likelihood}
p(Y_{t_{n+1}}|X_{t_{n+1}}) = \exp\Big(- \f{\big(Y_{t_{n+1}} - h(X_{t_{n+1}}) \big)^2}{R}\Big)
\end{equation}
is the likelihood function, and $p(Y_{t_{n+1}} | \mathcal{Y}_{t_n})$ in the denominator normalizes the filtering density at the time instant $t_{n+1}$. 

\vspace{1em}

Then, we carry out the above prediction-update procedure recursively to propagate the filtering density $p(X_{t}|Y_t)$ over time.

\section{Adaptive Kernel Approximation Approach}\label{Methodology}

In this paper, we solve the Fokker-Planck equation \eqref{Fokker-Planck} numerically to generate an approximation for the predicted filtering density $p(X_{t_{n+1}}|\mathcal{Y}_{t_{n}})$, and we apply the Bayesian inference \eqref{Bayes} to calculate the estimated (posterior) filtering density $p(X_{t_{n+1}}|\mathcal{Y}_{t_{n+1}})$. In what follows, we shall give detailed discussions on the computational framework that we construct to apply the adaptive kernel approximation method to solve the optimal filtering problem.

\subsection{Prediction through Fokker-Planck equation}
For convenience of presentation, we denote 
$$\mathcal{L}_{b, \sigma}p_t : =  - \sum_{i=1}^d \f{\p}{\p x_i}\big[ b_i(x, t) p(x, t)\big] + \sum_{i, j=1}^d \f{\p^2}{\p x_i \p x_j} \big[ D_{i, j} \ p(x, t)\big],$$
and we call $\mathcal{L}_{b, \sigma}$ the the Fokker-Planck operator in this paper. The Prediction step in our adaptive kernel approximation approach will focus on deriving a numerical solver for the Fokker-Planck equation
\begin{equation}\label{F-P}
\f{\p p(x, t)}{\p t} = \mathcal{L}_{b, \sigma}p,
\end{equation}
and the numerical solution to Eq. \eqref{F-P} will be our approximation to the predicted filtering density, which will be combined with the likelihood function to generate an estimated posterior filtering density.

Numerical methods for solving parabolic type PDEs, such like the Fokker-Planck equation, have been extensively studied \cite{Galerkin,Xiu-collocation, NSP3, Kushner-Dupuis}. However, when the dimension of the problem is high, solving Eq. \eqref{F-P} becomes an extremely expensive computational task \cite{Guannan-SG}. The primary challenge in obtaining numerical solutions to the Fokker-Planck equation is how to efficiently and effectively implement spatial dimensional approximation. Traditional mesh-based numerical methods, such like finite difference methods and finite element methods typically utilize polynomial approximations to describe solutions of the equation. However, due to the so-called ``curse of dimensionality'', the computational cost of solving the Fokker-Planck equation increases exponentially as the dimension of the problem increases.

In this work, we adopt the following kernel approximation scheme to approximate the solution of the Fokker-Planck equation
\begin{equation}\label{RBF}
p(x, t_n) \approx \sum_{k=1}^{K} \Phi^k_n(x), 
\end{equation}
where $K$ is the total number of kernels, and
\begin{equation}\label{Def:kernel}
\ds \Phi^k_n(x):=\omega^k_n \exp\Big(- \f{1}{2}\cdot (x - \mu^k_n)^{\top}(\Sigma^k_n)^{-1}(x - \mu^k_n)\Big)
\end{equation}
is a Gaussian type kernel function, which is parameterized by weight $\omega^k_n$, mean $\mu^k_n$ and covariance matrix $\Sigma^k_n$. Numerical analysis results have been derived to verify that the kernel approximation scheme \eqref{RBF} is capable of generating accurate approximations to wide range of function when the number of kernels $K$ is sufficiently large \cite{Bao_Meshfree, Kernel_Analysis_15, Kernel_Analysis_87}.  The reason that we pick Gaussian type kernels to approximate the target function $p$ is that $p$ is the filtering density, which describes a conditional probability distribution. In many situations in the optimal filtering problem, the filtering density is a bell-shaped function, which can be effectively approximated by Gaussian type functions. 

Then, assuming that we have a kernel approximation $p_{n}$ for the filtering density $p(X_{t_{n}}|\mathcal{Y}_{t_{n}})$, we introduce the following temporal discretization scheme to solve the Fokker-Planck equation
\begin{equation}\label{FP-Approx:Semi}
\tilde{p}_{n+1} = p_{n} + \mathcal{L}_{b, \sigma} p_{n} \cdot \Delta t_n,
\end{equation}
where $\Delta t_n : = t_{n+1} - t_n$ is the time step-size, and $\tilde{p}_{n+1}$ is a kernel approximation for the predicted filtering density $p(X_{t_{n+1}} \big| \mathcal{Y}_{t_n})$. Given kernels $\{\Phi^k_n\}_{k=1}^{K}$ for the approximated filtering density $\tilde{p}_{n}$ and the approximation scheme
\begin{equation}\label{Kernel-approximation}
p_{n}: =  \sum_{k=1}^{K} \Phi^k_n(x),
\end{equation}
we can rewrite Eq. \eqref{FP-Approx:Semi} as
\begin{equation}\label{FP-Approx:kernel}
\tilde{p}_{n+1} = \sum_{k=1}^{K} \Phi^k_n(x) + \Delta t_n \cdot \mathcal{L}_{b, \sigma} \Big( \sum_{k=1}^{K}\Phi^k_n(x) \Big),
\end{equation}
and we let
$$\tilde{p}_{n+1}:= \sum_{k=1}^{K} \tilde{\Phi}^k_{n+1},$$
where $\{\tilde{\Phi}^k_{n+1}\}_{k=1}^{K}$ is a set of kernels that approximates $\tilde{p}_{n+1}$.  We can see from the temporal discretization scheme \eqref{FP-Approx:kernel} that obtaining an approximation $\tilde{p}_{n+1}$ for the predicted filtering density $p(X_{t_{n+1}} \big| \mathcal{Y}_{t_n})$ is equivalent to finding parameters for kernels $\{\tilde{\Phi}^k_{n+1}\}_{k=1}^{K}$. Note that the kernels $\{ \Phi^k_{n}\}_{k=1}^{K}$ on the right hand side of Eq. \eqref{FP-Approx:kernel} are Gaussian (as introduced in Eq. \eqref{Def:kernel}). Hence the Fokker-Planck operator part, i.e. $\mathcal{L}_{b, \sigma} \Big( \sum_{k=1}^{K}\Phi^k_n(x) \Big)$ can be derived analytically. In this way, we transfer the computational cost of solving the Fokker-Planck equation from high dimensional spatial approximation to solving an optimization problem for kernel parameters. 

Since the target function for kernel approximation is a PDF, a relatively small number of Gaussian kernels may be sufficient to provide a reasonable description for the filtering density. On the other hand, solving the Fokker-Planck equation through Eq. \eqref{FP-Approx:kernel} suffers from the stability issue. When values of the drift function $b$ in the state equation Eq. \eqref{State} are large (or the time step-size $\Delta t_n$ is large), the drift term will generate a strong force that pushes the filtering density far from its current location. However, due to exponential decay of Gaussian tails, which would typically cause local behaviors of Gaussian kernels, the filtering density approximated by the kernel approximation scheme \eqref{Kernel-approximation} can only be transported to a limited distance. This can make our method difficult to track targets driven by  state equations with large drift terms.

\vspace{0.5em}

In the following subsection, we shall introduce an operator decomposition method to alleviate the above stability issue.

\vspace{0.5em}

\subsection{Operator decomposition}

The central idea of our operator decomposition method is to divide the Fokker-Planck operator into a drift operator and a diffusion operator. Then, we further decompose the drift operator into a linear component and a nonlinear component, and we provide analytical and numerical methods to characterize the linear component and the nonlinear component separately. 
\vspace{0.5em}

Before we introduce our decomposition strategy, we would like to point out the following facts of the Fokker-Planck operator:
\vspace{0.75em}

\noindent \textbf{Fact 1.} \textit{Given a PDF $p$, in the case that the diffusion coefficient $\sigma$ does not contain state $X$, we have
$$\mathcal{L}_{b, \sigma}p = \mathcal{L}_{b, 0} p + \mathcal{L}_{0, \sigma} p.$$
}

\vspace{0.75em}

\noindent \textbf{Fact 2.} \textit{The Fokker-Planck operator $\mathcal{L}_{b, \sigma}$ is linear, i.e., for two constants $a$, $b$ and two PDFs $p$, $q$, we have
$$\mathcal{L}_{b, \sigma}[a p + b q] = a \mathcal{L}_{b, \sigma} p + b \mathcal{L}_{b, \sigma} q.$$
Therefore, the kernel approximated filtering density under the Fokker-Planck operator can be written as
$$\mathcal{L}_{b, \sigma} p_{n} = \sum_{k=1}^{K} \mathcal{L}_{b, \sigma} \Phi^k_n(x),$$
and the right hand side of Eq. \eqref{FP-Approx:kernel} becomes
\begin{equation}\label{FP-Approx:kernel-linear}
\begin{aligned}
\sum_{k=1}^{L} \Phi^k_n(x) + \Delta t_n \cdot  \sum_{k=1}^{K}  \mathcal{L}_{b, \sigma} \Phi^k_n(x)
=  \sum_{k=1}^{K}\Big( \Phi^k_n(x) + \Delta t_n \cdot \mathcal{L}_{b, \sigma} \Phi^k_n(x) \Big) .
\end{aligned}
\end{equation}
The linearity of the Fokker-Planck operator allows us to discuss the propagation of each Gaussian kernel separately.}

\vspace{0.75em}

In light of \textbf{Fact 1}, we can handle the drift term first and then incorporate diffusion into the state propagation. \textbf{Fact 2} allows us to discuss state propagation kernel-by-kernel when necessary.
\vspace{0.75em}

In this work, instead of deriving the operator decomposition method directly under the numerical PDE framework, we first switch back to the state equation, and we consider the following Euler-Maruyama scheme that propagates each kernel $\Phi^k_n$ through the state equation
\begin{equation}\label{Discretize:SDE}
X^k_{n+1} = X^k_n + b(t_n, X^k_n) \Delta t_n + \sigma_{t_n} \Delta W_{t_n}, \quad k = 1, 2, \cdots, K,
\end{equation}
where the initial state $X^k_n \sim \Phi^k_n$, i.e. $X^k_n$ follows the distribution of the $k$-th Gaussian kernel, and $\Delta W_{t_n} := W_{t_{n+1}} - W_{t_n} \sim N(0, \Delta t_n \cdot \textbf{I}_d)$. In this way, by combining distributions for $\{X^k_{n+1} \}_{k=1}^{K}$ obtained through the discretized SDE scheme \eqref{Discretize:SDE}, we get a description for the predicted filtering density, which can also be considered as an approximation for the right hand side of Eq. \eqref{FP-Approx:kernel-linear}. 

To address the stability issue through operator decomposition and to transport Gaussian kernels effectively to the next time step, we introduce a linear approximation to the (nonlinear) drift function, and we denote it by $b^L(t_n, X^k_n) : = A X^k_n + \alpha$, where $A \in \mathbb{R}^{d \times d}$ and $\alpha \in \mathbb{R}^d$. The linear operator $b^L$ will be determined as the best linear approximation to $b$ in the sense of least square. In other words, we aim to find $A$ and $\alpha$ that will minimize the mean square error between the original drift function $b$ and the linear approximation $b^L$, i.e.
\begin{equation}\label{Optimization-linear}
\min_{A, \alpha} \E\left[ \Big( b(t_n, X^k_n) - (A X^k_n + \alpha) \Big)^2\right].
\end{equation}

To maintain the nonlinearity of the state dynamics, we introduce a residual function $b^N(t_n, X^k_n) := b(t_n, X^k_n) - b^L(t_n, X^k_n)$ that models the nonlinear component of $b$. Hence, the drift function is decomposed into a linear component $b^L$ and a nonlinear component $b^N$, i.e., $b(t_n, X^k_n) = b^L(t_n, X^k_n) + b^N(t_n, X^k_n)$, and the Euler-Maruyama scheme for the state equation can be interpreted as
\begin{equation*}
X^k_{n+1} = X^k_n + \big( b^L(t_n, X^k_n) + b^N(t_n, X^k_n)\big)\Delta t_n + \sigma \Delta W_{t_n}.
\end{equation*}

In what follows, we will introduce a three-step operator decomposition procedure to compute the predicted filtering density $\tilde{p}_{n+1}$.
\vspace{0.5em}

In the \textit{first step}, we only transport the filtering density via the linear component $b^{L}$. Specifically, we implement the following scheme
\begin{equation}\label{linear-transform}
\begin{aligned}
X^{k, L}_{n+1} := & X^k_n + b^L(t_n, X^k_n) \Delta t_n \\
= & (A \Delta t_n + \textbf{I} )X^k_n + \alpha \Delta t_n
\end{aligned}
\end{equation}
to propagate the filtering density at the time step $t_n$,
and we let 
$$T(X^k_n):=(A \Delta t_n + \textbf{I} )X^k_n + \alpha \Delta t_n$$
be the operator that formulates the linear component of the drift function, i.e. $X^{k, L}_{n+1}  = T(X^k_n)$. Note that a linear function will map a Gaussian distribution to a Gaussian distribution. 
Since $X^k_n$ follows a Gaussian distribution, $X^{k, L}_{n+1}$ will also follow a Gaussian distribution, which can be determined by the linear operator $T(\cdot)$, and we denote the distribution for $X^{k, L}_{n+1}$ by $p_{n+1}^{k, L}$.

\vspace{0.5em}

In the \textit{second step}, we incorporate the nonlinear component $b^N$ of the drift function to the filtering density so that both the linear and the nonlinear components are considered in the filtering density propagation.  Since $b^N$ does not linearly propagates $X^{k, L}_{n+1}$, we can not derive a Gaussian kernel directly from $p_{n+1}^{k, L}$ to obtain a kernel that describes the nonlinear component of the drift. 
In order to derive a kernel approximation for the predicted filtering density, which have considered the nonlinear component of the drift, we define an operator 
$$b^{N, T} (t_n, X^{k, L}_{n+1}):= b^N(t_n, T^{-1}(X^{k, L}_{n+1})) = b^N(t_n, X^k_n).$$ 
Then, with Gaussian distributions $\{p_{n+1}^{k, L}\}_{k=1}^{K}$ that describe random variables $\{X^{k, L}_{n+1}\}_{k=1}^K$ (introduced in Eq. \eqref{linear-transform}), we introduce the following PDE type solver to calculate a distribution $\hat{p}_{n+1}$ defined by
\begin{equation}\label{NL-FP}
\hat{p}_{n+1} = \sum_{k=1}^{K} \Big( p_{n+1}^{k, L} + \mathcal{L}_{b^{N, T} , 0} \ p_{n+1}^{k, L} \Delta t_n \Big),
\end{equation}
where $\mathcal{L}_{b^{N, T} , 0}$ is a Fokker-Planck operator with drift $b^{N, T}$, and the diffusion $\sigma$ is chosen as $0$. The PDF $\hat{p}_{n+1}$ on the left hand side of Eq. \eqref{NL-FP} is an approximation for the predicted filtering density before incorporation of the diffusion term, and we use kernel approximation scheme to represent $\hat{p}_{n+1}$, i.e. 
\begin{equation}\label{kernel:nonlinear}
\hat{p}_{n+1} = \sum_{k=1}^{K} \hat{\Phi}^k_{n+1}(x),
\end{equation}
where $\{\hat{\Phi}^k_{n+1}\}_{k=1}^{K}$ is a set of Gaussian kernels, and we will introduce the procedure to determine parameters for $\{\hat{\Phi}^k_{n+1}\}_{k=1}^{K}$ in the next subsection.

\vspace{0.5em}

Finally, in the \textit{third step} we add diffusion back to the predicted filtering density. 
Since we assume that the state dynamics are perturbed by additive noises in this work, for each Gaussian kernel $\hat{\Phi}^k_{n+1}$ that approximates $\hat{p}_{n+1}$ in Eq. \eqref{kernel:nonlinear}, we can simply introduce the extra diffusion information by adding $\Delta t_n \sigma_{t_n} \sigma_{t_n}^{\top}$ to the covariance of $\hat{\Phi}^k_{n+1}$ and get a kernel $\tilde{\Phi}^k_{n+1}$ to approximate the predicted filtering density at time stage $t_{n+1}$. As a result, we obtain the kernel approximation for the predicted filtering density as follows
\begin{equation}\label{Prior:approx}
\tilde{p}_{n+1} = \sum_{k=1}^{K} \tilde{\Phi}^k_{n+1}.
\end{equation}

\vspace{0.5em}

In the above three-step procedure, we can see that the first step and the third step can be implemented analytically, and the second step incorporates the nonlinear behavior of the dynamical model, which needs an optimization procedure to determine kernel parameters. In what follows, we will introduce an adaptive boosting algorithm to achieve this goal.


\subsection{Adaptive boosting algorithm for kernel training}
Recall that each kernel in Eq. \eqref{kernel:nonlinear} is Gaussian and has the expression 
$$\hat{\Phi}^{k}_{n+1}(x) = \hat{\omega}^k_{n+1} \exp\Big(- \f{1}{2}\cdot (x - \hat{\mu}^k_{n+1})^{\top}(\hat{\Sigma}^k_{n+1})^{-1}(x - \hat{\mu}^k_{n+1})\Big).$$
Our optimization procedure aims to find kernel parameters $\{( \hat{\omega}^k_{n+1},  \hat{\Sigma}^k_{n+1},  \hat{\mu}^k_{n+1})\}_{k=1}^{K}$ so that the left hand side of Eq. \eqref{NL-FP}, which is determined by $\{\hat{\Phi}^{k}_{n+1}\}_{k=1}^{K}$, will be equal to the right hand side, which is defined by the linear transformed Gaussian distributions $\{p^{k, L}_{n+1}\}_{k=1}^{L}$. We denote 
\begin{equation}\label{target}
g_{n+1} : = \sum_{k=1}^{K} \Big( p_{n+1}^{k, L} + \mathcal{L}_{b^{N, T} , 0}\ p_{n+1}^{k, L} \cdot \Delta t_n \Big)
\end{equation}
for convenience of presentation. Since $\{p_{n+1}^{k, L}\}_{k=1}^{K}$ are Gaussian functions, $g_{n+1}$ defined in Eq. \eqref{target} can be derived analytically.

\begin{table}[h!]\caption{Boosting Algorithm}\label{Algorithm:Boosting}
\vspace{0.5em}
\centering
\begin{tabular} {p{0.9\textwidth}}
\hline\noalign{\smallskip}
{\bf Algorithm 1}: {\em Boosting algorithm to adaptively generate kernels.}\\
\noalign
{\smallskip}\hline
\noalign{\smallskip}
\vspace{-0.1cm}
\begin{spacing}{0.8}
\begin{algorithmic}
\item[] Initialize the kernel approximation as $\hat{p}_{n+1}(x) = 0$; define target function $g_{n+1}$ through Eq. \eqref{target}; set global approximation tolerance $tol$.
\vspace{0.75em}
\item[\textbf{while}] $k =1, 2, \cdots, K$, \textbf{do} \\ \vspace{-0.5em}
\begin{description}
\item[-] \hspace{-0.25em} Generate $M$ global state samples, denoted by $\{\hat{x}_{n+1}^{(m)}\}_{m=1}^M$, from the kernel approximated distribution based on $\{p^{k, L}_{n+1}\}_{k=1}^{K}$. \vspace{-0.25em}
\item[-] \hspace{-0.25em}  Evaluate the approximation error on each state sample and calculate $e_m := g(\hat{x}_m) - \hat{p}_{n+1}(\hat{x}_m)$ for $m=1, 2, \cdots, M$. \vspace{-0.5em}
\item[-] \hspace{-0.25em}  Compute global error $\ds E_g = \frac{1}{M}\sum_{m=1}^M (e_m)^2$. If $E_g < tol$, break and set weights for other kernels $0$, i.e. $\omega_j = 0$, $k < j \leq K$. Otherwise, continue.
\item[-] \hspace{-0.25em} Locate the state sample with the largest approximation error, i.e. find $m^*$ s.t. $e_{m^*} = \max_m e_m$. \vspace{-0.25em}
\item[-] \hspace{-0.25em}  Generate a Gaussian kernel $\hat{\Phi}^{k}_{n+1}$ centered at the state sample that suffers from the largest error, i.e. choose the initial guess for the mean as $\hat{\mu}^k_{n+1} = \hat{x}_{m^*}$.\vspace{-0.25em}
\item[-] \hspace{-0.25em}  Solve a local optimization problem to determine the weight and covariance for the kernel  $\hat{\Phi}^{k}_{n+1}$ by comparing values of $\hat{\Phi}^{k}_{n+1}$ (treated as the left hand side of Eq. \eqref{NL-FP}) with $g_{n+1}$ on locally generated state samples near the kernel center $\hat{\mu}^k_{n+1}$. \vspace{-0.25em}
\item[-] \hspace{-0.25em} Add the locally trained kernel $\hat{\Phi}^{k}_{n+1}$ to kernel approximation $\hat{p}_{n+1}$, i.e. let $\hat{p}_{n+1}  = \hat{p}_{n+1} + \hat{\Phi}^{k}_{n+1}$. 
\end{description}
\item[\textbf{end while}]
\end{algorithmic}
\vspace{-1.2em}
\end{spacing}\\
\hline
\end{tabular}
\end{table}

In this work, instead of finding all the kernel parameters at the same time by solving a large scale optimization problem, we adopt the so-called ``boosting algorithm'', which sequentially minimizes the approximation error. Specifically, we introduce the Boosting Algorithm in Table \ref{Algorithm:Boosting} to determine the parameter set $\{( \hat{\omega}^{k}_{n+1},  \hat{\Sigma}^{k}_{n+1},  \hat{\mu}^{k}_{n+1})\}_{k=1}^{K}$.
\vspace{1em}

The boosting algorithm introduced in Table \ref{Algorithm:Boosting} will adaptively generate kernels, and this adaptive kernel approximation procedure allows us to capture more important features (modes) in the filtering density. Also, the Gaussian tails of the kernels can provide reasonable description for low density regions in the filtering density, which will make our method stable. 

%

\subsection{Bayesian update for filtering density}

To incorporate the observational information to the predicted filtering density, we apply Bayesian inference \eqref{Bayes}. Since the predicted filtering density is described by multiple kernels, we apply Bayesian inference to each Gaussian kernel and obtain a kernel for the posterior filtering density. Specifically, for each state point $x$, let 
$$ p^{k, \text{post}}_{n+1}(x) = \tilde{\Phi}^k_{n+1}(x) p(Y_{t_{n+1}} \big| x),$$
where $p(Y_{t_{n+1}}\big| x)$ is the likelihood function introduced in Eq. \eqref{Likelihood} with a given state position $X_{t_{n+1}} = x$ and $\tilde{\Phi}$ is a Gaussian kernel in Eq. \eqref{Prior:approx} that approximates the predicted filtering density $\tilde{p}_{n+1}$. In this way, the entire posterior filtering density is approximated by
\begin{equation}\label{posterior}
p^{\text{post}}_{n+1} = \sum_{k=1}^K p^{k, \text{post}}_{n+1}.
\end{equation}

Note that each kernel $p^{k, \text{post}}_{n+1}$ that we use to approximate the overall posterior filtering density $p^{\text{post}}_{n+1}$ may not be Gaussian due to the nonlinear observation. To derive an approximation by Gaussian kernels, we train a new set of Gaussian kernels to describe the posterior filtering density. Specifically, we introduce a kernel approximation 
$$p_{n+1} := \sum_{k=1}^{K} \Phi^{k}_{n+1},$$ 
and we let $p_{n+1}$ be an approximation to the approximated posterior filtering density $p^{\text{post}}_{n+1}$, i.e. $p_{n+1} \approx p^{\text{post}}_{n+1}$. To this end, we adopt the same Boosting Algorithm framework  introduced in Table \ref{Algorithm:Boosting} again to adaptively generate Gaussian kernels $\{\Phi^{k}_{n+1}\}_{k=1}^{K}$, and a normalization procedure will be implemented to $\{\Phi^{k}_{n+1}\}_{k=1}^{K}$ to make $p_{n+1}$ a PDF.

\subsection{Summary of the algorithm}

In this subsection, we summarize our algorithm in Table \ref{Algorithm}.

\begin{table}[h!]\caption{Summary of the algorithm}\label{Algorithm}
\vspace{0.5em}
\centering
\begin{tabular} {p{0.9\textwidth}}
\hline\noalign{\smallskip}
{\bf Algorithm 2}: {\em Algorithm of the adaptive kernel method.}\\
\noalign
{\smallskip}\hline
\noalign{\smallskip}
\vspace{-0.1cm}
\begin{spacing}{0.8}
\begin{algorithmic}
\item[] Initialize the filtering density $p_{0}$ with kernels $\{\Phi^k_0\}_{k=1}^K$.
\vspace{0.75em}
\item[\textbf{For}] $n =0, 1, 2, 3, \cdots $  \\
 \begin{description}
 \item[$\bullet$] Prediction Step:
\begin{description}
\item[-] \hspace{-0.25em} Generate $\{p^{k, L}_{n+1}\}_{k=1}^{K}$ through Eq. \eqref{linear-transform} to incorporate the linear component (determined through Eq. \eqref{Optimization-linear} ) of the drift function. \vspace{-0.25em}
\item[-] \hspace{-0.25em} Use the Boosting Algorithm described in Table \ref{Algorithm:Boosting} to incorporate the nonlinear component of the drift function and generate Gaussian kernels $\{\hat{\Phi}^k_{n+1}\}_{k=1}^K$ from $g_{n+1}$ (defined in Eq. \eqref{target} ) to approximate $\hat{p}_{n+1}$. \vspace{-0.25em}
\item[-] \hspace{-0.25em} Add $\Delta t_n \sigma_{t_n} \sigma_{t_n}^{\top}$ to the covariance of each Gaussian kernel $\hat{\Phi}^k_{n+1}$ to incorporate state diffusion and get the kernel $\tilde{\Phi}^k_{n+1}$ to approximate the predicted filtering density $\tilde{p}_{n+1}$ via Eq. \eqref{Prior:approx}. 
\end{description}
\item[$\bullet$] Update Step
\begin{description}
\item[-] \hspace{-0.25em} Carry out Bayesian inference to generate a posterior filtering density $p^{\text{post}}_{n+1}$ defined in Eq. \eqref{posterior}. \vspace{-0.25em}
\item[-] \hspace{-0.25em} Carry out Boosting Algorithm in Table \ref{Algorithm:Boosting} again to obtain a Gaussian kernel approximation $p_{n+1} = \{\Phi_{n+1}^k\}_{k=1}^{K}$ to approxiomate $p^{\text{post}}_{n+1}$. \vspace{-0.25em}
\item[-] \hspace{-0.25em} Normalize $\{\Phi_{n+1}^k\}_{k=1}^{K}$ to make $p_{n+1}$ a PDF, and $p_{n+1}$ is the estimated filtering density at time stage $n+1$. \vspace{-0.25em}
\end{description}
\end{description}
\item[\textbf{end}]
\end{algorithmic}
\vspace{-1.2em}
\end{spacing}\\
\hline
\end{tabular}
\end{table}

\section{Numerical Experiments}\label{Numerics}

In this section, we present three numerical examples to demonstrate the performance of our adaptive kernel method for solving the optimal filtering problem. We first present a demonstration example to show how our adaptive kernel approximation method will adaptively capture the main features of the filtering density in state propagation. In the second example, we solve a benchmark optimal filtering problem, i.e. the bearing-only tracking problem, and we compare our method with the particle filter method \cite{APF} and the ensemble Kalman filter method \cite{Evense_EnKF} to show accuracy and efficiency of the adaptive kernel method.  Then, in Example 3 we solve a high dimensional Lorenz-96 tracking problem, which is a well-known challenging optimal filtering problem due to the chaotic behavior of the state model.

\subsection{Example 1: Demonstration for adaptive kernel approximation.}
We use the first numerical example to demonstrate the performance of our adaptive kernel approximation method in propagating state dynamics. Instead of solving an entire optimal filtering problem, we only present the effectiveness of our method in transporting a probability distribution through the Fokker-Planck equation, and the primary computational effort of our approach lies on using kernels to approximate the Fokker-Planck operator. Since the filtering density is approximated by Gaussian kernels, Gaussian type diffusions can be directly added to the target distribution. Therefore, in this example we shall focus on the drift part of the Fokker-Planck operator, i.e. $\mathcal{L}_{b, 0}$, and the drift term is defined by the following 2D function: 
\begin{equation*}
    b(x_1,x_2)  = 
    \begin{bmatrix}
    x_2^2\\
    0
    \end{bmatrix}
    +
    \begin{bmatrix}
    3&4\\
    3&2
    \end{bmatrix}
    \begin{bmatrix}
    x_1\\
    x_2
    \end{bmatrix}
    +
    \begin{bmatrix}
    3\\
    -2
    \end{bmatrix}.
\end{equation*}

\begin{figure}[h!]
    \centering
    \includegraphics[width=0.9\textwidth]{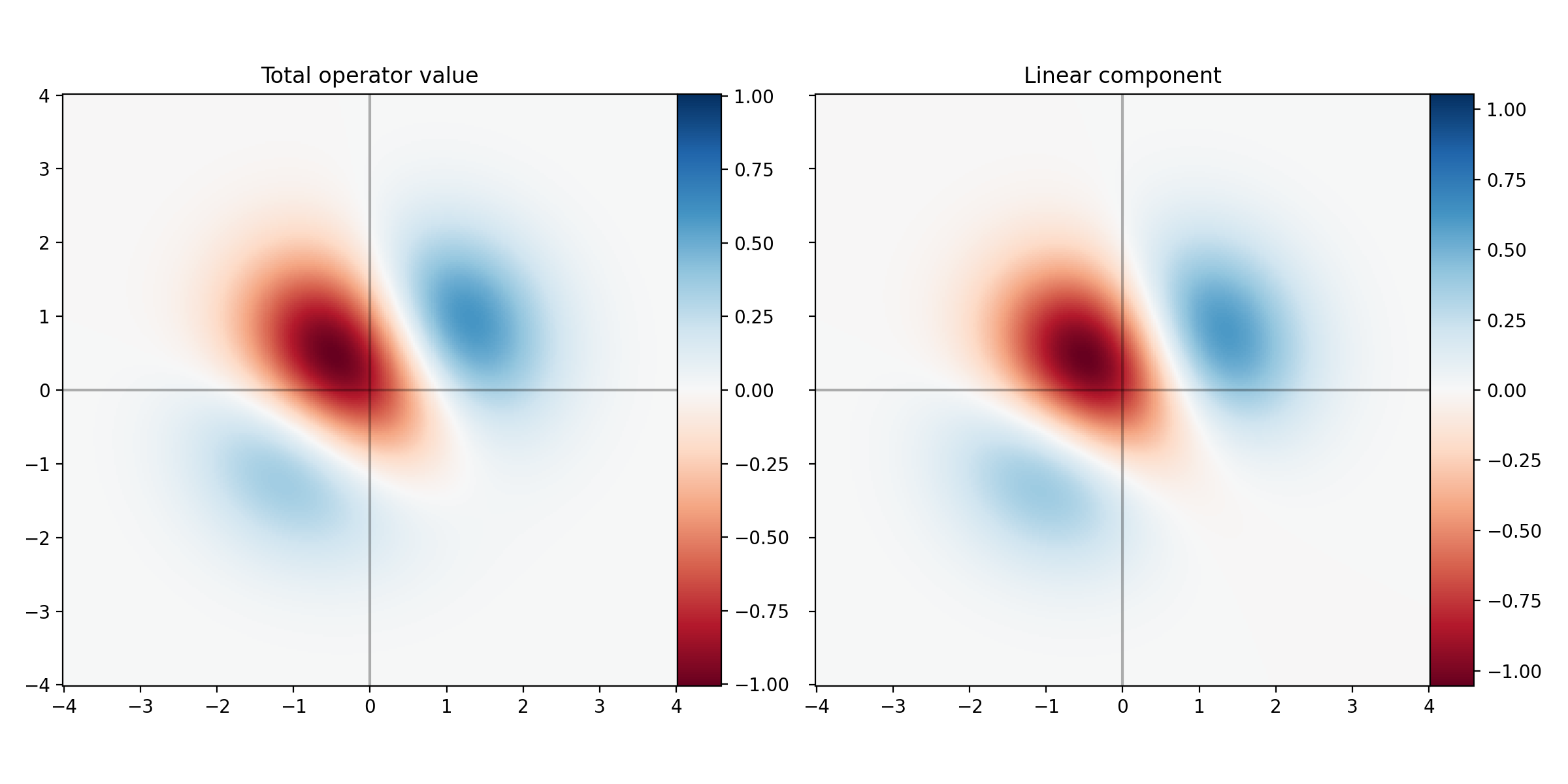}
    \caption{Example 1. Linear component in describing the Fokker-Planck operator}
    \label{Ex1: decomposition}
\end{figure}

For convenience of presentation, we consider state propagation in time with step-size be $1$. Then, we choose the initial distribution as a standard Gaussian distribution, denoted by $\Phi$, and we apply the drift operator $\mathcal{L}_{b, 0}$ to $\Phi$. In this way, the target function that we try to use our kernel method to approximate is $F := \Phi + \mathcal{L}_{b, 0} \Phi$.  In Figure \ref{Ex1: decomposition}, we present the original target function $F$ driven by the operator $\mathcal{L}_{b, 0} $ on left, and the linear approximation for function $F$ obtained by the linear transportation Eq. \eqref{linear-transform} is presented on the right. From this figure, we can see that the linear component can roughly capture the main feature of the target function. 

\begin{figure}[h!]
    \centering
    \includegraphics[width=0.9\textwidth]{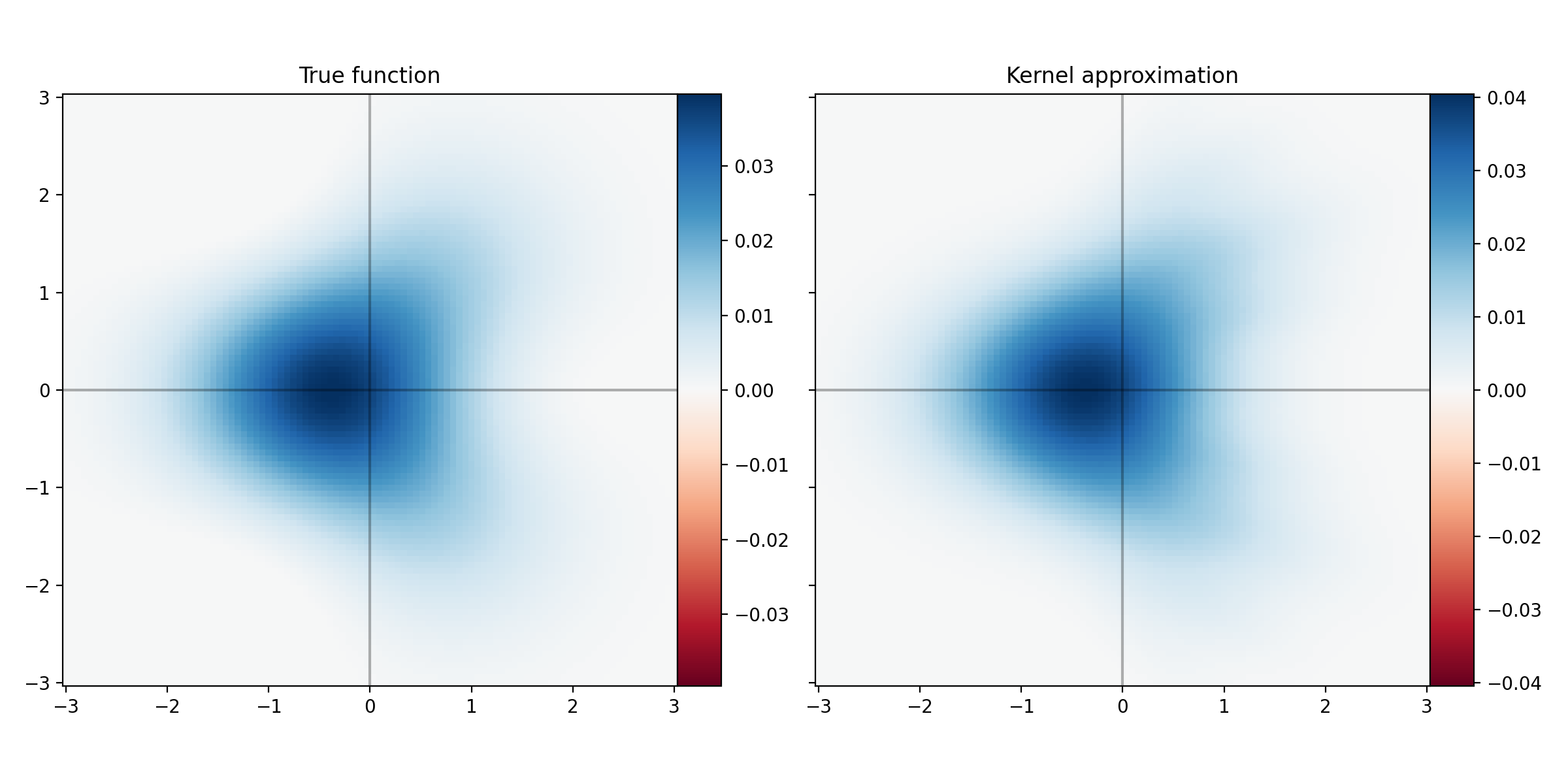}
    \caption{Example 1. Accuracy of approximation obtained by the boosting algorithm.}
    \label{Ex1: Kernel_fitting}
\end{figure}

To demonstrate the performance of kernel method in approximating the nonlinear component (described in Eq. \eqref{NL-FP}) of the operator and the effectiveness of the adaptive boosting algorithm, we compare the analytically derived true nonlinear component of the function with the approximated nonlinear component in Figure \ref{Ex1: Kernel_fitting}.  
The subplot on the left shows the true function that we aim to approximate, and the subplot on the right is our approximated function by using the boosting algorithm introduced in Table \ref{Algorithm:Boosting}. We use blue-to-red colors to represent function values, and we can see from this figure that the boosting  algorithm can accurately capture the true function, which describes the nonlinear component of the Fokker-Planck operator.

To show more details of the performance of the adaptive kernel construction in the boosting algorithm, 
\begin{figure}[h!]
    \centering
    \includegraphics[width=0.9\textwidth]{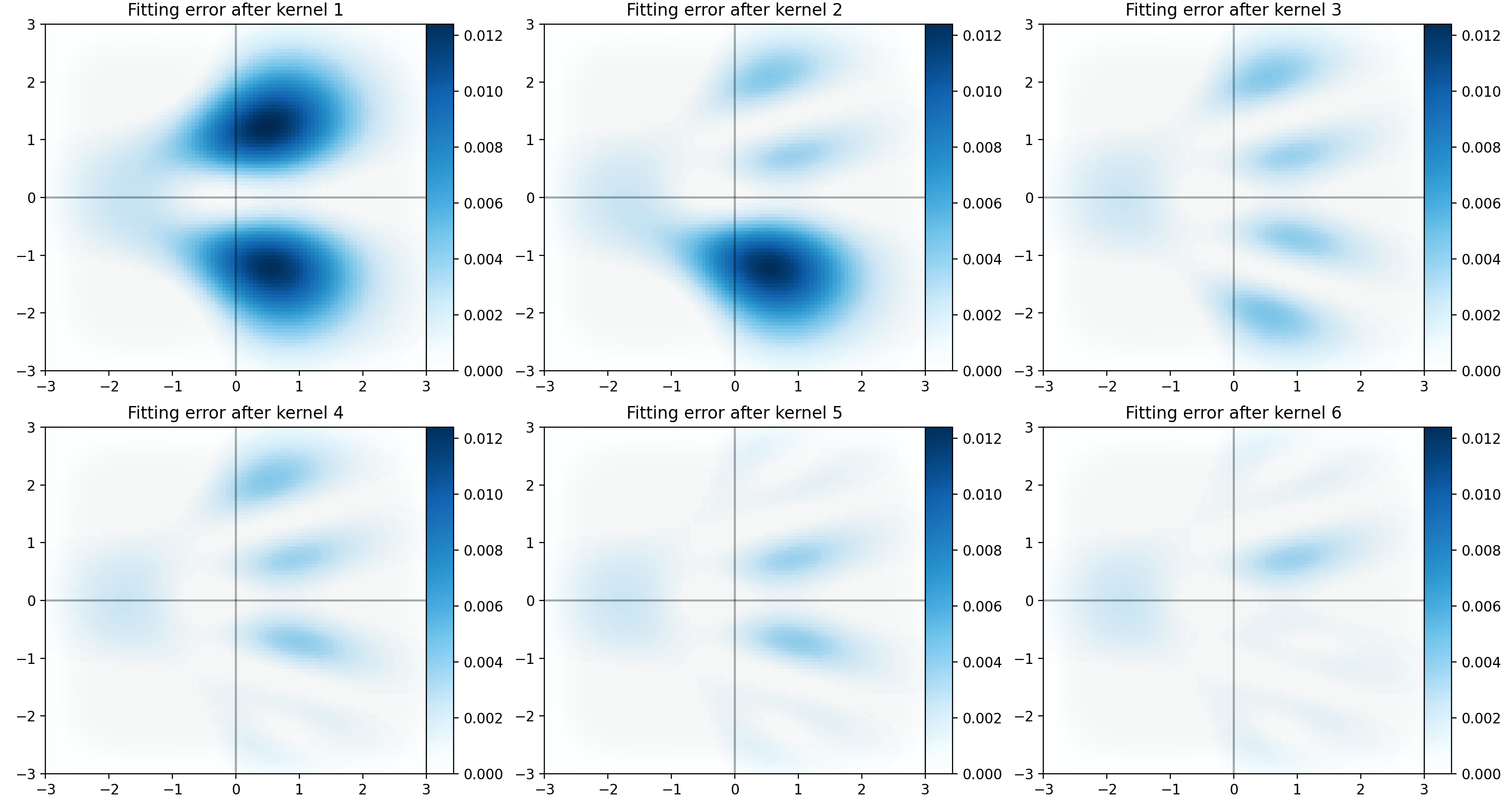}
    \caption{Example 1. Performance of the adaptive boosting algorithm in reducing approximation errors in fitting the nonlinear component of the operator.}
    \label{Ex1: Adaptive_Kernel}
\end{figure}
we present the approximation errors after fitting up to $6$ kernels in Figure \ref{Ex1: Adaptive_Kernel}. From this figure, we can see that by using only one kernel to describe the nonlinear component of the Fokker-Planck operator, the main part of the function in the region $[-1, 0] \times [-1, 1]$ (presented in the left subplot in Figure \ref{Ex1: Kernel_fitting}) is well fitted, and two remaining features that represent two tails in the function (plotted in Figure \ref{Ex1: Kernel_fitting}) need to be fitted. Then, by adding the second and the third kernels, we can successfully approximate those two tails and get low overall fitting errors. As more and more kernels are added, we get rid of higher error regions one-by-one. As a result, we obtain more and more accurate approximations to the nonlinear component of the drift operator.

\subsection{Example 2: Bearing-only tracking}
In this example, we solve the bearing-only tracking problem, which is a benchmark optimal filtering problem in practice. Specifically, we aim to track a moving target driven by the following state dynamics
\begin{equation}\label{Ex2:State}
    dX_t  = \begin{bmatrix}
    v^1_t\\
    v^2_t\\
    0\\
    0
    \end{bmatrix}dt
    + 
    \begin{bmatrix}
    \sigma_1 & 0 & 0 & 0 \\
    0 & \sigma_2 & 0 & 0 \\
    0 & 0 & \sigma_3  & 0 \\
    0 & 0 & 0 & \sigma_4 
    \end{bmatrix}
    \begin{bmatrix}
    dW^1_t\\
    dW^2_t\\
    dW^3_t\\
    dW^4_t
    \end{bmatrix},
\end{equation}
where $X_t = [x^1_t, x^2_t, v^1_t, v^2_t]^{\top}$, $[x_t^1, x_t^2]^{\top}$ describes the 2D location of the target, and $v_t^1$, $v^2_t$ are the velocities in $x_1$ and $x_2$ directions, respectively. $W_t = [W_t^1, W_t^2, W_t^3, W_t^4]$ is a 4D Brownian motion that brings uncertainty to the state model, which is driven by the diffusion coefficient $\sigma :=  \begin{bmatrix}
    \sigma_1 & 0 & 0 & 0 \\
    0 & \sigma_2 & 0 & 0 \\
    0 & 0 & \sigma_3  & 0 \\
    0 & 0 & 0 & \sigma_4 
    \end{bmatrix}$.
\vspace{0.5em}

In order to estimate the location of the target, we place two detectors on different observation platforms to collect bearing angles as observational data. Specifically, the observational data is given by the following observational function
\vspace{0.25em}
\begin{equation}\label{Ex2:Observation}
Y_t^i = \arctan\big( \f{x_t^2 - x^2_{\text{i-platform}}}{x_t^1 - x^1_{\text{i-platform}}}\big) + \xi_i, \qquad i = 1, 2,
\end{equation}
where $(x^1_{\text{i-platform}}, x^2_{\text{i-platform}})^{\top}$ gives the location of the $i$-th platform, and $\xi_i$ is the observation noise of the $i$-th detector.  
\begin{figure}[h!]
    \centering
    \includegraphics[width=0.7\textwidth]{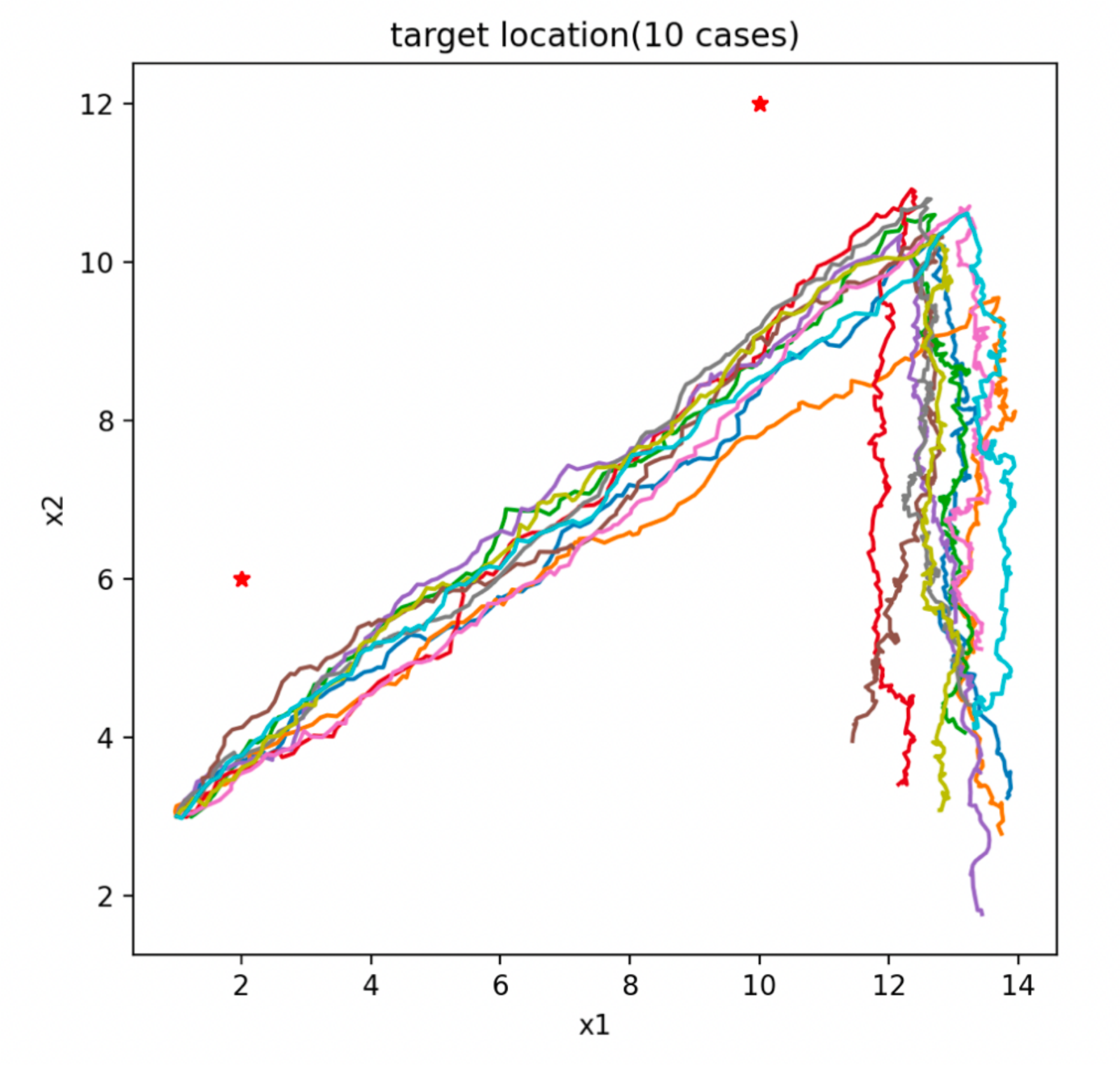}
    \caption{Example 2. Demonstration of $10$ sample trajectories of the target.}
    \label{Ex2: Sample_paths}
\end{figure}

In this example, we track the target over the time period $t \in [0, 3]$ with initial state $X_0 = [1, 3, 10, 6]^{\top}$, and we let $\Delta t = 0.01$, i.e. we track $300$ time steps. The diffusion coefficient is chosen as $\sigma_1 = \sigma_2 = 0.5$, $\sigma_3 = \sigma_4 = 0.3$, and we locate two platforms at $(2, 6)^{\top}$ and $(10, 12)^{\top}$, respectively.  To demonstrate the stability of our method compared with other state-of-the-art methods, we assume that there's an unexpected turn in the target moving direction at the time instant $t = 1.2$, which would challenge the robustness of optimal filtering methods. In Figure \ref{Ex2: Sample_paths}, we plot $10$ sample target trajectories (by using different colors), and we mark the observation platforms with red stars. From this figure, we can see that the target is designed to move in front of the observation platforms and then it makes a sharp turn downward. 

\begin{figure}[h!]
    \centering
    \includegraphics[width=0.8\textwidth]{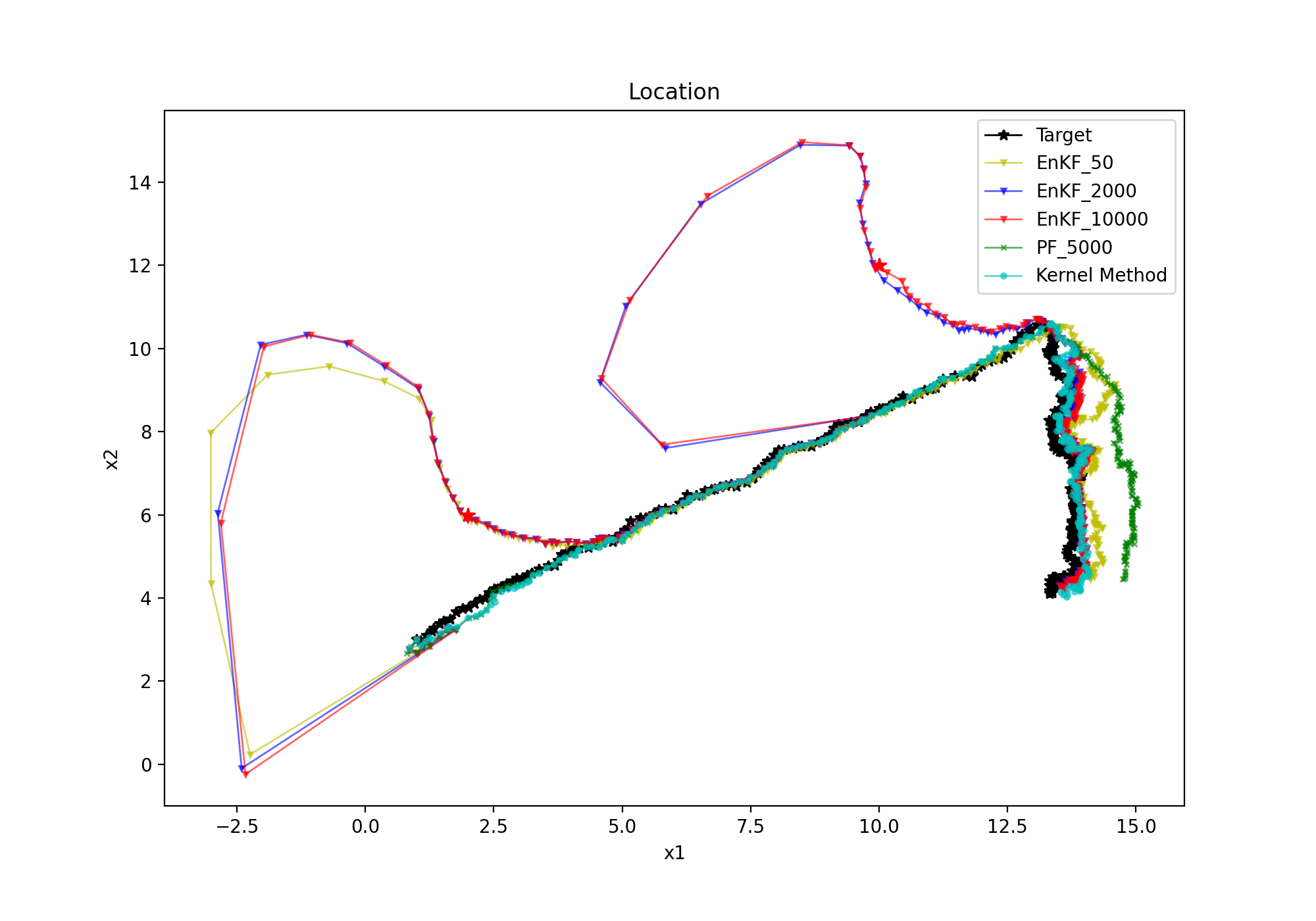}
    \caption{Example 2. Comparison of tracking performance in solving the bearing-only tracking problem.}
    \label{Ex2: Tracking_Comparison}
\end{figure}

In Figure \ref{Ex2: Tracking_Comparison}, we present a comparison experiment, in which we compare the tracking accuracy between our adaptive kernel method with two state-of-the-art optimal filtering methods, i.e. the ensemble Kalman filter and the particle filter.
To implement our adaptive kernel method, we use up to $20$ kernels to approximate the filtering density, and active kernels are adaptively selected by the boosting algorithm (described in Table \ref{Algorithm:Boosting}). For the ensemble Kalman filter, we choose $50$, $2,000$ and $10,000$ realizations of Kalman filter samples to implement this tracking task. In the particle filter, we use $5000$ particles to generate empirical distributions for the filtering density.  In the figure, we use the black curve (marked by stars) to represent a sample of real target trajectory and use other colored curves to represent the estimates obtained by various optimal filtering methods. The yellow, blue, and red curves (marked by triangles) are estimates for the target location obtained by using the ensemble Kalman filter (EnKF) with $50$, $2,000$, and $10,000$ realizations of Kalman filter samples, respectively. The green curve (marked by crosses) gives the particle filter (PF) estimates (obtained by using $5,000$ particles). The cyan curve (marked by dots) describes the estimates obtained by our adaptive kernel method. 

From this figure, we can see that the EnKF doesn't provide accurate estimates for the target location when the target is right below a detector -- no matter how many realizations of samples we use in the EnKF. The poor performance of the EnKF is caused by the high nonlinearity of observational data (bearing angles introduced in Eq. \eqref{Ex2:Observation}) when the target moves in front of detectors. For the PF, we can see that it provides accurate estimates until the sharp turn at the time instant $t =1.2$. Then, the PF loses track of the target due to the degeneracy of particles when trying to adjust the change of the target location. On the other hand, the kernel method always keeps on track, and it gives accurate estimates all the time during the tracking period.

\begin{figure}[h!]
    \centering
    \includegraphics[width=0.8\textwidth]{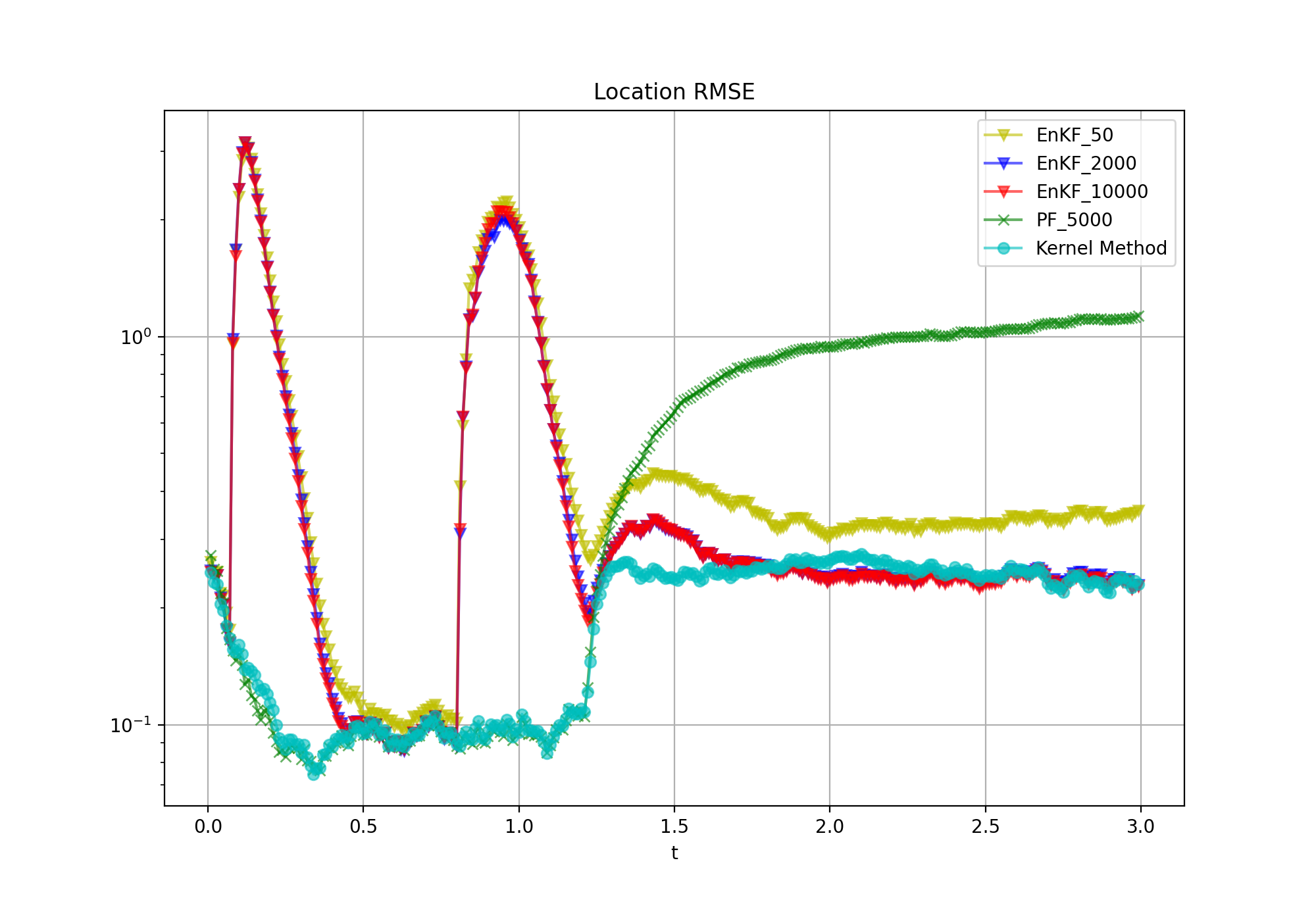}
    \caption{Example 2. Comparison of root mean square errors (RMSEs) with respect to time.}
    \label{Ex2: RMSE}
\end{figure}

To further examine the performance of different optimal filtering methods in solving the bearing-only tracking problem \eqref{Ex2:State}-\eqref{Ex2:Observation}, we repeat the above experiment $100$ times and calculate the root mean square errors (RMSEs) of target tracking performance. The $\log$ scaled RMSEs of each method with respect to time are presented in Figure \ref{Ex2: RMSE}. From this figure, we can see that the adaptive kernel method (cyan curve marked by dots) has the lowest RMSEs, and it can provide good accuracy even after the sharp turn of the target. The PF (green curve marked by  crosses) has low RMSEs at first. However, the errors increase dramatically at the turning point of the target trajectories. On the other hand, the EnKF estimates (yellow, blue and red curves marked by triangles) always suffer from low accuracy when the target passes the detectors. But the EnKF can recover quickly from inaccurate estimates, which indicates that the EnKF is a more stable method compared with the PF.

To summarize the general performance of each method, we present the accumulated RMSEs (the combined RMSEs over the tracking period) together with the  CPU time of each method (average over the above $100$ repeated tests) in Table \ref{CPU_Comparison}. 
\renewcommand{\arraystretch}{1.25}
\begin{table} [h!]\small
\leftmargin=6pc \caption{Performance comparison} \label{CPU_Comparison} 
\begin{center}
\begin{tabular}{|c|c|c|c|c|c|}
 \hline   &EnKF 50&  EnKF 2,000   & EnKF 10,000 & PF 5000 & Kernel Method \\
\hline   Accumulated RMSEs & $158.99$ & $134.75$ & $134.87$ & $169.58$ & $56.75$ \\
\hline   CPU time (seconds) & $0.32$ & $11$ & $56$ & $70$  & $50$\\
\hline
\end{tabular}\end{center}
\end{table}
The CPU that we use is a AMD Ryzen 5 5600X with 6 core 12 processing threads. We can see from this table that the PF has the lowest accuracy with the highest computational cost, which is mainly caused by the degeneracy of particles. The EnKF can solve the problem with very low computational cost. However, the accuracy of the EnKF does not improve much even we use a lot more realizations of Kalman filter samples.   The kernel learning method, on the other hand,  has much lower RMSEs compared with the EnKF and the PF with moderate cost.

\subsection{Example 3: Lorenz-96 model}
To examine the performance of the adaptive kernel method in solving high dimensional problems, in this example we solve the Lorenz-96 tracking problem, which is a benchmark high dimensional optimal filtering problem. The state model is given by the following  stochastic dynamical system
\begin{equation}\label{Lorenz-96}
    x_t^i = ((x_t^{i+1} - x_t^{i-2})x_t^{i-1} - x_t^{i} + F)dt + \sigma^i dW_t^i, \quad i = 1, 2, \cdots, d
\end{equation}
where $X_t = [x_t^1, x_t^2 \cdots, x_t^d]^{\top}$ is the target state. In the Lorenz-96 model \eqref{Lorenz-96}, we let $x_t^{-1}=x_t^{d-1}$, $x_t^{0}=x_t^{d}$, $x_t^{1}=x_t^{d+1}$, $W_t = \{W_t^1, W_t^2 \cdots, W_t^d\}$ is a $d$-dimensional Brownian motion, and $\sigma = [\sigma^1, \sigma^2 \cdots, \sigma^d]^{\top}$ is the diffusion coefficient. It is well-known that when $F = 8$, the Lorenz-96 model has chaotic behavior, which makes the corresponding optimal filtering problem very challenging. In this example, we track the state $X$ of the Lorenz-96 model over the time period $t \in [0, 3]$, and we let $d = 10$. As a commonly used scenario when tracking the Lorenz-96 model, we simulate the Lorenz-96 model with time step-size $\Delta t = 0.001$, and we assume that we receive data of the state with time step-size $\Delta t = 0.1$. Therefore, the Bayesian inference procedure is implemented after every $100$ simulation steps. In other words, we carry out one update step in every $100$ predication steps.

\begin{figure}[h!]
   \hspace{-4em} \includegraphics[width=1.2\textwidth]{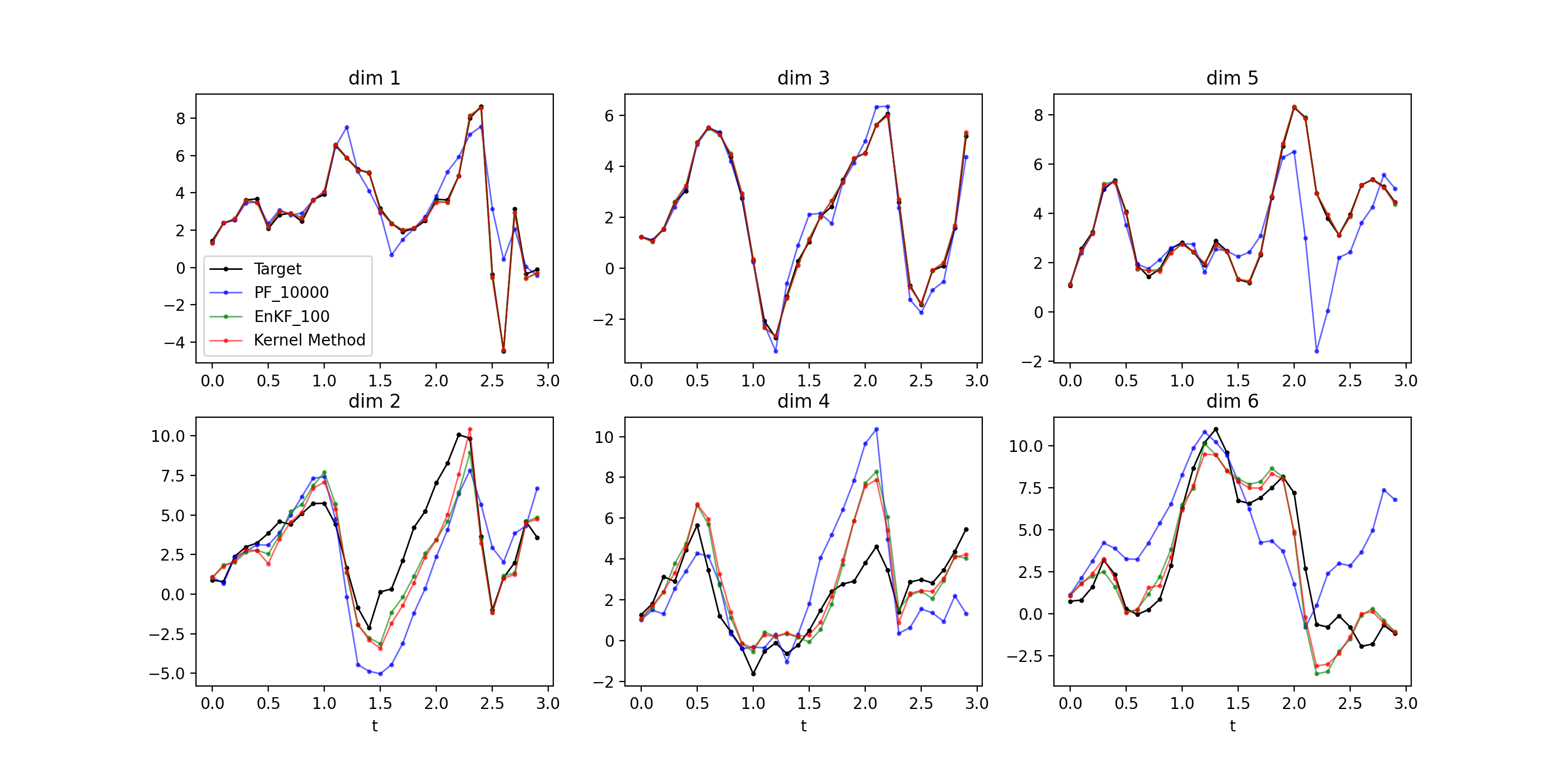}
    \caption{Example 3. Comparison of tracking performance in dimensions $1$ to $6$.}
    \label{Ex3: Comparison}
\end{figure}

In this example, the observational data that we receive to estimate the state of the Lorenz-96 model are noise perturbed direct state observations in odd dimensions, i.e. 
$$Y_{t} = [x_t^1, x_t^3, x_t^5, x_t^7, x_t^9]^{\top} + \xi_t,$$
where $\xi_t \sim N(0, 1)$ is a standard Gaussian noise.

In Figure \ref{Ex3: Comparison}, we compare our kernel method with the PF and the EnKF, and we present the state estimation performance of each method in the first six dimensions. In each subplot, the black curve shows a sample of real target trajectory of the Lorenz-96 model state. The blue curve is the PF estimates obtained by using $10,000$ particles to represent the empirical distribution of the state. The green curve is the EnKF estimates obtained by using $100$ realizations of Kalman filter samples. The red curve gives the estimates obtained by our kernel method, and we use at most $20$ kernels to approximate the filtering density in the adaptive boosting algorithm when fitting the nonlinear component of the state drift. From this figure, we can see that the EnKF has comparable estimation performance to the kernel method due to the linear observations, and the usage of ``ensemble estimation'' in the EnKF can handle the nonlinearity of the state dynamics. On the other hand, the PF provides low tracking accuracy. Especially, the long simulation period (without an update) in this example would cause more severe degeneracy issue since no data can be used to resample the particles.

\begin{figure}[h!]
    \centering
    \includegraphics[width=0.8\textwidth]{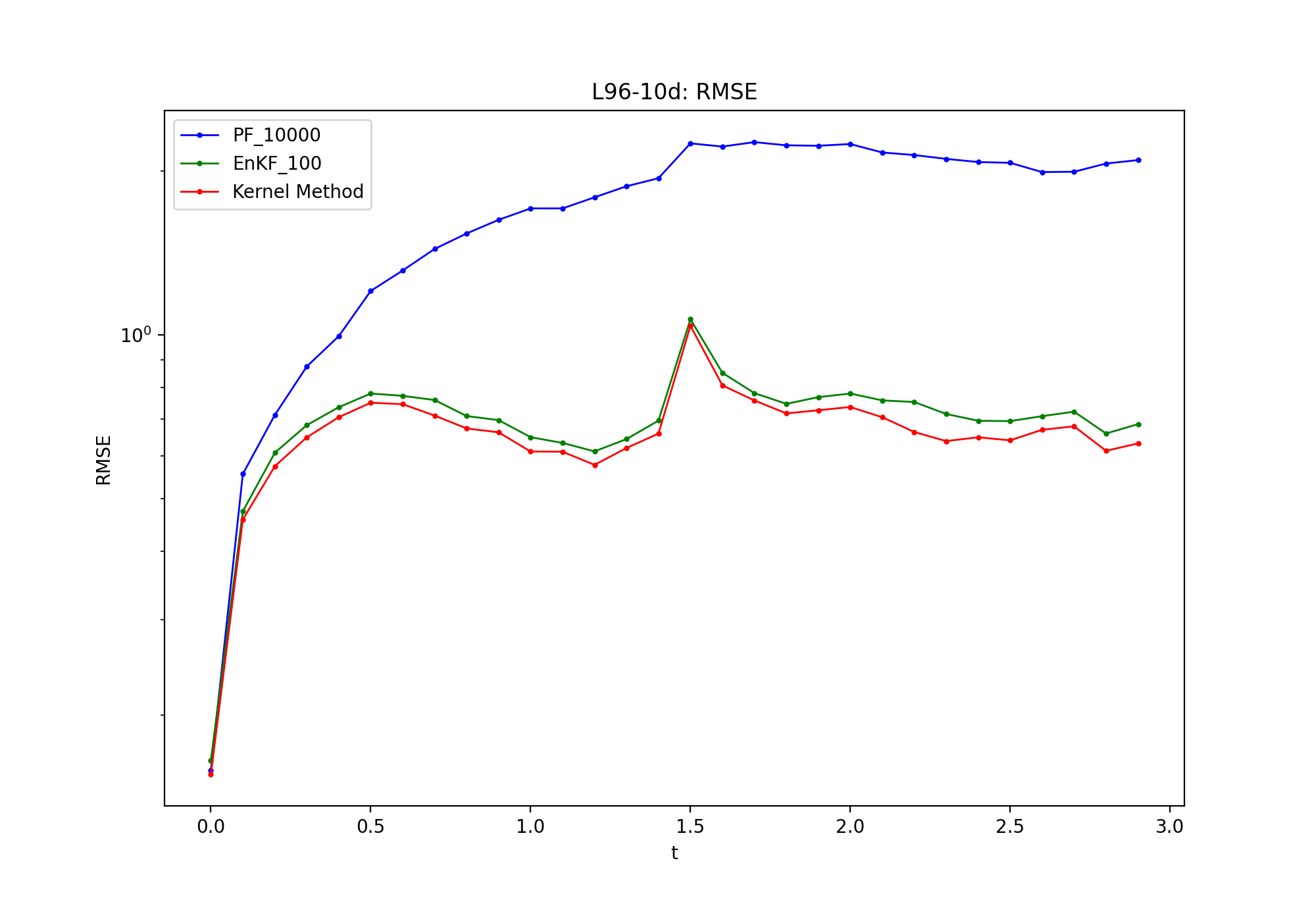}
    \caption{Example 3. Comparison of RMSEs with respect to time.}
    \label{Ex3: RMSE}
\end{figure}

To confirm the comparison result presented in Figure \ref{Ex3: Comparison}, we repeat the above experiment $100$ times and present the $\log$ scaled RMSEs of each method with respect to time in Figure \ref{Ex3: RMSE}. We can see from this figure that the PF has much higher errors compared with the EnKF and the kernel method while both the EnKF and the kernel method have similar RMSEs in this Lorenz-96 tracking problem. 

\section{Summary and conclusions}\label{Conclusion}
In this paper, we developed an adaptive kernel method to solve the optimal filtering problem. The main idea of our method is to use a set of Gaussian kernels to approximate the filtering density of a target dynamical state model.  Due to the fact that the filtering density describes a probabilistic distribution, Gaussian kernels can effectively characterize the distribution, which is often a bell-shaped function. Then, an operator decomposition method is introduced to efficiently propagate the state of the model, and adaptive boosting algorithm is applied to adaptively capture important features of the filtering density.

Three numerical experiments are presented to examine the performance of our kernel method. In the first example, we presented the effectiveness of the adaptive kernel method in characterizing propagation of the filtering density. In the second example and the third example, we compared the performance of the kernel method with two state-of-the-art methods, i.e. the particle filter and the ensemble Kalman filter, in solving benchmark optimal filtering problems. Results in our numerical experiments indicate that our method has high accuracy and high stability advantage compared with the particle filter, and it outperforms the ensemble Kalman filter when data provide highly nonlinear state observations.

\section*{Acknowledgement}

This work is partially supported by U.S. Department of Energy through FASTMath Institute and Office of Science, Advanced Scientific Computing Research program under the grant DE-SC0022297. The second author (FB) would also like to acknowledge the support from U.S. National Science Foundation through project DMS-2142672.

\bibliographystyle{plain}


\begin{thebibliography}{10}

\bibitem{MCMC-PF}
C.~Andrieu, A.~Doucet, and R.~Holenstein.
\newblock Particle markov chain monte carlo methods.
\newblock {\em J. R. Statist. Soc. B}, 72(3):269--342, 2010.

\bibitem{Bao_Kernel_Arxiv}
R.~Archibald and F.~Bao.
\newblock A kernel learning method for backward sde filter.
\newblock {\em arXiv: 2201.10600}, 1, 2021.

\bibitem{Bao_Kernel}
R.~Archibald and F.~Bao.
\newblock Kernel learning backward sde filter for data assimilation.
\newblock {\em Journal of Computational Physics}, 455:111009, 2022.

\bibitem{Bao_filter_jump}
F.~Bao, Y.~Cao, and H.~Chi.
\newblock Adjoint forward backward stochastic differential equations driven by
  jump diffusion processes and its application to nonlinear filtering problems.
\newblock {\em Int. J. Uncertain. Quantif.}, 9(2):143--159, 2019.

\bibitem{Bao_CMS}
F.~Bao, Y.~Cao, and X.~Han.
\newblock Forward backward doubly stochastic differential equations and optimal
  filtering of diffusion processes.
\newblock {\em Communications in Mathematical Sciences}, 18(3):635--661, 2020.

\bibitem{Bao_first}
F.~Bao, Y.~Cao, A.~J. Meir, and W.~Zhao.
\newblock A first order scheme for backward doubly stochastic differential
  equations.
\newblock {\em SIAM/ASA J. Uncertain. Quantif.}, 4(1):413--445, 2016.

\bibitem{Bao_Zakai}
F.~Bao, Y.~Cao, C.~Webster, and G.~Zhang.
\newblock A hybrid sparse-grid approach for nonlinear filtering problems based
  on adaptive-domain of the {Z}akai equation approximations.
\newblock {\em SIAM/ASA J. Uncertain. Quantif.}, 2(1):784--804, 2014.

\bibitem{Bao_Half}
F.~Bao, Y.~Cao, and W.~Zhao.
\newblock Numerical solutions for forward backward doubly stochastic
  differential equations and zakai equations.
\newblock {\em International Journal for Uncertainty Quantification},
  1(4):351--367, 2011.

\bibitem{Bao_Semi}
F.~Bao, Y.~Cao, and W.~Zhao.
\newblock A first order semi-discrete algorithm for backward doubly stochastic
  differential equations.
\newblock {\em Discrete and Continuous Dynamical Systems-Series B}, 5(2):1297
  -- 1313, 2015.

\bibitem{BCZ_2018}
F.~Bao, Y.~Cao, and W.~Zhao.
\newblock A backward doubly stochastic differential equation approach for
  nonlinear filtering problems.
\newblock {\em Commun. Comput. Phys.}, 23(5):1573--1601, 2018.

\bibitem{BSDE_filter}
F.~Bao and V.~Maroulas.
\newblock Adaptive meshfree backward {SDE} filter.
\newblock {\em SIAM J. Sci. Comput.}, 39(6):A2664--A2683, 2017.

\bibitem{Bao_CiCP_2019}
Feng Bao, Richard Archibald, and Petro Maksymovych.
\newblock Backward {SDE} filter for jump diffusion processes and its
  applications in material sciences.
\newblock {\em Communications in Computational Physics}, 27:589--618, 2020.

\bibitem{NSP3}
A.~Davie and J.~Gaines.
\newblock Convergence of numerical schemes for the solution of the parabolic
  stochastic partial differential equations.
\newblock {\em Math. Comp.}, 70:121--134, 2001.

\bibitem{BSDE_finance}
N.~El~Karoui, S.~Peng, and M.~C. Quenez.
\newblock Backward stochastic differential equations in finance.
\newblock {\em Math. Finance}, 7(1):1--71, 1997.

\bibitem{EvensenBook}
G.~Evensen.
\newblock {\em Data assimilation: the ensemble {K}alman filter}.
\newblock Springer, 2006.

\bibitem{Evense_EnKF}
G.~Evensen.
\newblock The ensemble {K}alman filter for combined state and parameter
  estimation: {M}onte {C}arlo techniques for data assimilation in large
  systems.
\newblock {\em IEEE Control Syst. Mag.}, 29(3):83--104, 2009.

\bibitem{Gobet-Zakai}
Emmanuel Gobet, Gilles Pag{\`e}s, Huy{\^e}n Pham, and Jacques Printems.
\newblock Discretization and simulation of the {Z}akai equation.
\newblock {\em SIAM J. Numer. Anal.}, 44(6):2505--2538 (electronic), 2006.

\bibitem{Multi-Kernel_Learning}
Mehmet G\"{o}nen and Ethem Alpayd\i~n.
\newblock Multiple kernel learning algorithms.
\newblock {\em J. Mach. Learn. Res.}, 12:2211--2268, 2011.

\bibitem{particle-filter}
N.J Gordon, D.J Salmond, and A.F.M. Smith.
\newblock Novel approach to nonlinear/non-gaussian bayesian state estimation.
\newblock {\em IEE PROCEEDING-F}, 140(2):107--113, 1993.

\bibitem{Galerkin}
W.~Grecksch and P.~E. Kloeden.
\newblock Time-discretised {G}alerkin approximations of parabolic stochastic
  {PDE}s.
\newblock {\em Bull. Austral. Math. Soc.}, 54(1):79--85, 1996.

\bibitem{Kernel_learning}
Thomas Hofmann, Bernhard Sch\"{o}lkopf, and Alexander~J. Smola.
\newblock Kernel methods in machine learning.
\newblock {\em Ann. Statist.}, 36(3):1171--1220, 2008.

\bibitem{Bao_Meshfree}
Sun Hui and Feng Bao.
\newblock Meshfree approximation for stochastic optimal control problems.
\newblock {\em Communications in Mathematical Research}, 37(3):387--420, 2021.

\bibitem{Kalman1961}
R.~E. Kalman and R.~S. Bucy.
\newblock New results in linear filtering and prediction theory.
\newblock {\em Transactions of the ASME--Journal of Basic Engineering},
  83(Series D):95--108, 1961.

\bibitem{Kang-PF}
K.~Kang, V.~Maroulas, I.~Schizas, and F.~Bao.
\newblock Improved distributed particle filters for tracking in a wireless
  sensor network.
\newblock {\em Comput. Statist. Data Anal.}, 117:90--108, 2018.

\bibitem{Kushner-Dupuis}
H.J. Kushner and P.~Dupuis.
\newblock Numerical methods for stochastic control problems in continuous time.
\newblock In {\em Applications of Mathematics}, volume~24. Springer-Verlag, New
  York, 1992.

\bibitem{NIPS1999_Boosting}
Llew Mason, Jonathan Baxter, Peter Bartlett, and Marcus Frean.
\newblock Boosting algorithms as gradient descent.
\newblock In S.~Solla, T.~Leen, and K.~M\"{u}ller, editors, {\em Advances in
  Neural Information Processing Systems}, volume~12. MIT Press, 2000.

\bibitem{PP1994}
{\'E}tienne Pardoux and Shi~Ge Peng.
\newblock Backward doubly stochastic differential equations and systems of
  quasilinear {SPDE}s.
\newblock {\em Probab. Theory Related Fields}, 98(2):209--227, 1994.

\bibitem{APF}
Michael~K. Pitt and Neil Shephard.
\newblock Filtering via simulation: auxiliary particle filters.
\newblock {\em J. Amer. Statist. Assoc.}, 94(446):590--599, 1999.

\bibitem{Kernel_Analysis_87}
M.~J.~D. Powell.
\newblock {\em Radial Basis Functions for Multivariable Interpolation: A
  Review}, page 143–167.
\newblock Clarendon Press, USA, 1987.

\bibitem{Sny-particle}
C.~Snyder, T.~Bengtsson, P.~Bickel, and J.~Anderson.
\newblock Obstacles to high-dimensional particle filtering.
\newblock {\em Mon. Wea. Rev.}, 136:4629--4640, 2008.

\bibitem{Kernel_Analysis_15}
Dougal~J. Sutherland and Jeff Schneider.
\newblock On the error of random fourier features.
\newblock UAI'15, Arlington, Virginia, USA, 2015. AUAI Press.

\bibitem{vanLeeuwen}
P.~J. van Leeuwen.
\newblock Nonlinear data assimilation in geosciences: an extremely efficient
  particle filter.
\newblock {\em Q. J. Roy. Meteor. Soc.}, 136(653):1991--1999, 2010.

\bibitem{DA-applications}
Bin Wang, Xiaolei Zou, and Jiang Zhu.
\newblock Data assimilation and its applications.
\newblock {\em Proceedings of the National Academy of Sciences},
  97(21):11143--11144, 2000.

\bibitem{Xiu-collocation}
Dongbin Xiu and Jan~S. Hesthaven.
\newblock High-order collocation methods for differential equations with random
  inputs.
\newblock {\em SIAM J. Sci. Comput.}, 27(3):1118--1139 (electronic), 2005.

\bibitem{zakai}
Moshe Zakai.
\newblock On the optimal filtering of diffusion processes.
\newblock {\em Z. Wahrscheinlichkeitstheorie und Verw. Gebiete}, 11:230--243,
  1969.

\bibitem{Guannan-SG}
Guannan Zhang and Max Gunzburger.
\newblock Error analysis of a stochastic collocation method for parabolic
  partial differential equations with random input data.
\newblock {\em SIAM J. Numer. Anal.}, 50(4):1922--1940, 2012.

\bibitem{Zhang_Zakai}
H.~Zhang and D.~Laneuville.
\newblock Grid based solution of zakai equation with adaptive local refinement
  for bearing-only tracking.
\newblock {\em IEEE Aerospace Conference}, 2008.

\end{thebibliography}

\end{document}